# OPTIMAL AND EFFICIENT CROSSOVER DESIGNS FOR COMPARING TEST TREATMENTS WITH A CONTROL TREATMENT

By A. S. Hedayat[1] and Min Yang[2]

*University of Illinois at Chicago and University of Nebraska-Lincoln*

This paper deals exclusively with crossover designs for the purpose of comparing $t$ test treatments with a control treatment when the number of periods is no larger than $t+1$. Among other results it specifies sufficient conditions for a crossover design to be simultaneously A-optimal and MV-optimal in a very large and appealing class of crossover designs. It is expected that these optimal designs are highly efficient in the entire class of crossover designs. Some computationally useful tools are given and used to build assorted small optimal and efficient crossover designs. The model robustness of these newly discovered crossover designs is discussed.

**1. Introduction.** Crossover designs, where experimental units are used in two or more ($p$) periods for the purpose of evaluating and studying two or more ($t$) treatments, have proven effective in a wide range of applications in agriculture, sensory testing and especially in phase I and phase II pharmaceutical clinical trials. The rigorous study of these designs and their optimality and efficiency has a history of at least 30 years. These optimality and efficiency results are almost all for those situations where all treatment comparisons are equally important. Some published works include Hedayat and Afsarinejad (1975, 1978), Cheng and Wu (1980), Kunert (1983, 1984),

Received November 2002; revised April 2004.
[1]Supported by NSF Grant DMS-01-03727, National Cancer Institute Grant RO1-CA099904 and NIH Grant P50-AT00155 (jointly supported by the National Center for Complementary and Alternative Medicine, the Office of Dietary Supplements, the Office for Research on Women's Health and the National Institute of General Medicine). The contents are solely the responsibility of the authors and do not necessarily represent the official views of NIH.
[2]Supported by NSF Grant DMS-03-04661.

Hedayat and Zhao (1990), Stufken (1991), Carrière and Reinsel (1993), Matthews (1994), Kushner (1997, 1998), Kunert and Stufken (2002) and Hedayat and Stufken (2003). Recently Hedayat and Yang (2003, 2004) have obtained additional results for balanced uniform designs and the designs suggested by Stufken (1991). We refer the readers to the excellent expository review paper by Stufken (1996) for additional references.

In many pharmaceutical studies, experimenters are more interested in the comparisons between $t$ experimental test treatments and an established standard or control treatment. In our study we shall designate the class of all such designs based on $t$ test treatments, a control treatment and $n$ experimental units each used in $p$ periods by $\Omega_{t+1,n,p}$. Unfortunately, the published literature contains very few useable results related to this very important problem. Perhaps the associated mathematical difficulties could be a major reason for this lack of useful published results. Majumdar (1988) considered A-optimal and MV-optimal crossover designs for comparing several test treatments with a control treatment and established that some known strongly balanced crossover designs can be combinatorially modified to obtain optimal designs for this problem. But Majumdar's result is limited to the situation where $t = w^2$, $p = m(w^2 + w)$ and $m \geq 2$. The first nontrivial case will be $t = 4$ and $p = 12$. In crossover designs, for a variety of reasons, usually it is undesirable to have large numbers of periods. In many cases, the experimenter is interested in designs in which the number of periods is less than or equal to the number of total treatments, that is, $p \leq t + 1$. When $p = 2$, Hedayat and Zhao (1990) studied this problem and obtained useful results.

Throughout this article, unless otherwise specified, we always assume $3 \leq p \leq t + 1$. The main purpose of this paper is to identify and construct crossover designs which are simultaneously A-optimal and MV-optimal in some desirable class of designs in terms of $t$ test treatments versus a control treatment. Section 2 introduces the model and the notation. Section 3 contains four lemmas which provide useful mathematical tools for research on optimal crossover designs in general. The main result is presented in Section 4. Some assorted examples are given in Section 5. Section 6 discusses further issues and indicates some important open problems. The Appendix contains a series of propositions which are used in Sections 4 and 5.

**2. Model of response.** Selection of an appropriate model of response for the data gathered under a crossover design in $\Omega_{t+1,n,p}$ is rather difficult and complicated. While several models have been introduced in the literature, the model which has been mostly entertained by design theorists is the traditional homoscedastic, additive and fixed effects model which, in the notation of Hedayat and Afsarinejad (1975), can be written as

$$Y_{dku} = \mu + \alpha_k + \beta_u + \tau_{d(k,u)} + \rho_{d(k-1,u)} + e_{ku},$$

(2.1)
$$k = 1, \dots, p;\ u = 1, \dots, n,$$

where $Y_{dku}$ denotes the response from unit $u$ in period $k$ to which treatment $d(k,u)$ was assigned. In this model $\mu$ is the general mean, $\alpha_k$ is the effect due to period $k$, $\beta_u$ is the effect due to unit $u$, $\tau_{d(k,u)}$ is the direct effect of $d(k,u)$, $\rho_{d(k-1,u)}$ is the carryover or residual effect of treatment $d(k-1,u)$ on the response observed on unit $u$ in period $k$ (by convention $\rho_{d(0,u)} = 0$) and the $e_{ku}$'s are independently normally distributed errors with mean 0 and variance $\sigma^2$. Recently another model has been introduced into the literature by Afsarinejad and Hedayat (2002) and some optimal designs have been identified. Kunert and Stufken (2002) used this new model of Afsarinejad and Hedayat and obtained additional optimal designs under this model.

Hereafter we shall designate the $t$ test treatments by $1, 2, \dots, t$ and the control treatment by 0. Throughout this paper, for each design $d$, we adopt the notation $n_{diu}$, $\tilde{n}_{diu}$, $l_{dik}$, $m_{dij}$, $r_{di}$, $\tilde{r}_{di}$ and $\hat{r}_{d0}$ to denote the number of times that treatment $i$ is assigned to unit $u$, the number of times this happens in the first $p-1$ periods associated with $u$, the number of times treatment $i$ is assigned to period $k$, the number of times treatment $i$ is immediately preceded by treatment $j$, the total replication of treatment $i$ in the $n$ experimental units, the total replication of treatment $i$ limited to the first $p-1$ periods in the design and the total replication of control treatment 0 limited to periods 2 to $p$, respectively.

In matrix notation we can write model (2.1) for the $np$ total observations as

$$Y_d = \mu 1 + P\alpha + U\beta + T_d\tau_d + F_d\rho_d + e, \qquad (2.2)$$

where $Y_d = (Y_{d11}, Y_{d21}, \dots, Y_{dpn})'$, $\alpha = (\alpha_1, \dots, \alpha_p)'$, $\beta = (\beta_1, \dots, \beta_n)'$, $\tau_d = (\tau_0, \dots, \tau_t)'$, $\rho_d = (\rho_0, \dots, \rho_t)'$, $e = (e_{11}, e_{21}, \dots, e_{pn})'$, $P = \mathbb{1}_n \otimes I_p$, $U = I_n \otimes \mathbb{1}_p$, $T_d = (T'_{d1}, \dots, T'_{dn})'$ and $F_d = (F'_{d1}, \dots, F'_{dn})'$. Here $T_{du}$ stands for the $p \times (t+1)$ period-treatment incidence matrix for unit $u$ under design $d$ and $F_{du} = LT_{du}$ with the $p \times p$ matrix $L$ defined as

$$\begin{pmatrix} 0_{1\times(p-1)} & 0 \\ I_{(p-1)\times(p-1)} & 0_{(p-1)\times 1} \end{pmatrix}.$$

The information matrix $C_d$, for direct effects $\tau_d = (\tau_0, \dots, \tau_t)'$, can now be expressed as

$$C_d = T'_d \operatorname{pr}^{\perp}([P|U|F_d])T_d,$$

where $\operatorname{pr}^{\perp}(X) = I - \operatorname{pr}(X)$ and $\operatorname{pr}(X) = X(X'X)^{-}X'$.

**3. Preliminary lemmas.** For comparing test treatments with a control, the most frequently used optimality criterion is A-optimality. An A-optimal design minimizes $\sum_{i=1}^{t} \mathrm{Var}_d(\hat{\tau}_i - \hat{\tau}_0)$. Let $M_d$ be the information matrix for the contrast vector $(\tau_1 - \tau_0, \ldots, \tau_t - \tau_0)'$. Then an A-optimal design for this contrast vector is a design which minimizes $\mathrm{Tr}(M_d^{-1})$. Another optimality criterion, which is associated with A-optimality, is MV-optimality. An MV-optimal design minimizes $\mathrm{Max}_{i=1,\ldots,t} \mathrm{Var}_d(\hat{\tau}_i - \hat{\tau}_0)$. The following well-known lemma indicates that an A-optimal design is also an MV-optimal design.

LEMMA 1. *An A-optimal design is also an MV-optimal design if its information matrix, $M_d$, is a completely symmetric matrix.*

Lemma 2 points out a well-known algebraic relation between the two information matrices $C_d$ and $M_d$.

LEMMA 2. *The information matrix $M_d$ can be obtained from the information matrix $C_d$ by*

$$M_d = V' C_d V, \tag{3.1}$$

*where $V$ is the $(t+1) \times t$ matrix, defined as*

$$\begin{pmatrix} 0_{1\times t} \\ I_{t\times t} \end{pmatrix}.$$

*Thus, $M_d$ can be simply obtained from $C_d$ by deleting the first row and the first column of $C_d$.*

From Lemma 2, we can explore the algebraic properties of $M_d$ via $C_d$. But unfortunately, $C_d$ for our problem is very complicated to deal with. It is known that when all test treatments and the control treatment are uniform in periods, that is, test treatments and the control treatment appear equally often in each period, we can ignore the impact of the period effects in $C_d$. This property is highly desirable when all test treatments and the control treatment are equally important. But when only the comparisons between test treatments and the control treatment are considered, this property is no longer desirable. In the latter case it is expected that the control treatment should have more replications than each test treatment. Lemma 3 shows that we can still ignore the period effects even if the replication of each test treatment and the control treatment are different. We state this result without a proof.

LEMMA 3. *For any crossover design $d$,*

$$C_d \le T_d' \mathrm{pr}^{\perp}(U) T_d - T_d' \mathrm{pr}^{\perp}(U) F_d (F_d' \mathrm{pr}^{\perp}(U) F_d)^{-} F_d' \mathrm{pr}^{\perp}(U) T_d, \tag{3.2}$$

*with equality for a crossover design with $l_{dik} = r_{di}/p$, $i = 0, \ldots, t$.*

Notice that the uniformity in periods is just a special case of the stated condition.

To find an A-optimal design, one needs to find a design that minimizes $\mathrm{Tr}(M_d^{-1})$. One way to achieve this is first to find the $\min_{d\in\Omega_{t+1,n,p}} \mathrm{Tr}(M_d^{-1})$ and then characterize the design that achieves this minimum value. Although Lemma 2 can help in simplifying the calculation of $M_d^{-1}$, it is still difficult to find the $\min_{d\in\Omega_{t+1,n,p}} \mathrm{Tr}(M_d^{-1})$. The main difficulty is that for a general design $d$, one cannot write down the explicit expression of $\mathrm{Tr}(M_d^{-1})$. Lemma 4 shows a method for finding the explicit expression of the achievable lower bound for $\mathrm{Tr}(M_d^{-1})$ for a very broad and general class of crossover designs.

LEMMA 4. *For any design $d$ in which $0<\tilde r_{d0}<n(p-1)$, we have*

$$\mathrm{Tr}(M_d^{-1}) \ge \frac{t-1}{x_0}+\frac{1}{y_0},$$

*where*

$$\begin{aligned}
x_0 = {} & \frac{t(np-r_{d0}-(1/p)\sum_{u=1}^n\sum_{i=1}^t n_{diu}^2)-(r_{d0}-(1/p)\sum_{u=1}^n n_{d0u}^2)}{t(t-1)} \\
& -\Bigg\{tp\Bigg(\sum_{i=1}^t\Bigg(m_{dii}-\frac{1}{p}\sum_{u=1}^n n_{diu}\tilde n_{diu}\Bigg) \\
& \qquad +\frac{1}{t}\Bigg(\hat r_{d0}-\frac{p-1}{p}r_{d0}-m_{d00}+\frac{1}{p}\sum_{u=1}^n n_{d0u}\tilde n_{d0u}\Bigg)\Bigg)^2\Bigg\} \\
& \times\Bigg\{(t-1)\Bigg[n(p-1)(pt-t-1)-(pt-t+p-2)\tilde r_{d0}+\sum_{u=1}^n \tilde n_{d0u}^2\Bigg]\Bigg\}^{-1}, \\
y_0 = {} & \frac{1}{t}\Bigg(r_{d0}-\frac{1}{p}\sum_{u=1}^n n_{d0u}^2\Bigg) \\
& -\Bigg\{p[n(p-1)-\tilde r_{d0}]\Bigg(m_{d00}-\frac{1}{p}\sum_{u=1}^n n_{d0u}\tilde n_{d0u}\Bigg)^2 \\
& \qquad +p\tilde r_{d0}\Bigg(\hat r_{d0}-\frac{p-1}{p}r_{d0}-m_{d00}+\frac{1}{p}\sum_{u=1}^n n_{d0u}\tilde n_{d0u}\Bigg)^2\Bigg\} \\
& \times\Bigg\{t\Bigg[np(p-1)\tilde r_{d0}-\tilde r_{d0}^2-n(p-1)\sum_{u=1}^n \tilde n_{d0u}^2\Bigg]\Bigg\}^{-1}.
\end{aligned}$$

*Further,* $\mathrm{Tr}(M_d^{-1})=\frac{t-1}{x_0}+\frac{1}{y_0}$ *will hold when the following conditions are satisfied:*

(i) $l_{dik} = r_{di}/p$, $i = 0, \dots, t$,

(ii) $T_d' \operatorname{pr}^\perp(U) T_d$, $T_d' \operatorname{pr}^\perp(U) F_d$ *and* $F_d' \operatorname{pr}^\perp(U) F_d$ *are invariant under any permutation of test treatments,*

(iii) *each test treatment appears at most once in each of $n$ first $p-1$ periods.*

PROOF. It can be shown that the elements in $T_d' \operatorname{pr}^\perp(U) T_d$ at positions $(i,i)$ and $(i,j)$ $(i \neq j)$ are $r_{di} - \frac{1}{p}\sum_{u=1}^n n_{diu}^2$ and $-\frac{1}{p}\sum_{u=1}^n n_{diu} n_{dju}$, respectively; the element in $T_d' \operatorname{pr}^\perp(U) F_d$ at position $(i,j)$ is $m_{dij} - \frac{1}{p}\sum_{u=1}^n n_{diu}\tilde{n}_{dju}$; the elements in $F_d' \operatorname{pr}^\perp(U) F_d$ at positions $(i,i)$ and $(i,j)$ $(i \neq j)$ are $\tilde{r}_{di} - \frac{1}{p}\sum_{u=1}^n \tilde{n}_{diu}^2$ and $-\frac{1}{p}\sum_{u=1}^n \tilde{n}_{diu}\tilde{n}_{dju}$, respectively.

Let $S_1 = I_t, S_2, \dots, S_N$ be the set of all $t \times t$ permutation matrices on the set of $t$ test treatments, where $N = t!$. Define $\bar{M}_d = \frac{1}{N}\sum_{i=1}^N S_i' M_d S_i$. By Lemma 2.2 of Majumdar and Notz (1983), $\bar{M}_d$ is superior to $M_d$ under the A-criterion. Define

$$\widetilde{S_i} = \begin{pmatrix} 1 & 0_{1\times t} \\ 0_{t\times 1} & S_i \end{pmatrix}, \qquad i = 1, \dots, N.$$

Then $\widetilde{S_i} = I_{t+1}, \widetilde{S}_2, \dots, \widetilde{S}_N$ will constitute the set of all $(t+1)\times(t+1)$ permutation matrices on the set of all treatments leaving the control treatment unchanged. It can be easily verified that $S_i' V' C_d V S_i = V' \widetilde{S_i}' C_d \widetilde{S_i} V$.

By Lemma 2, we have

$$\bar{M}_d = \frac{1}{N}\sum_{i=1}^N S_i' V' C_d V S_i = \frac{1}{N}\sum_{i=1}^N V' \widetilde{S_i}' C_d \widetilde{S_i} V = V' \left\{ \frac{1}{N}\sum_{i=1}^N \widetilde{S_i}' C_d \widetilde{S_i} \right\} V. \tag{3.3}$$

By (3.2) and Proposition 1 of Kunert and Martin (2000), we have

$$\begin{aligned} \sum_{i=1}^N \widetilde{S_i}' C_d \widetilde{S_i} &\le \sum_{i=1}^N \widetilde{S_i}' (T_d' \operatorname{pr}^\perp(U) T_d \\ &\qquad - T_d' \operatorname{pr}^\perp(U) F_d (F_d' \operatorname{pr}^\perp(U) F_d)^- F_d' \operatorname{pr}^\perp(U) T_d) \widetilde{S_i} \\ &\le \sum_{i=1}^N \widetilde{S_i}' (T_d' \operatorname{pr}^\perp(U) T_d) \widetilde{S_i} \\ &\quad - \left( \sum_{i=1}^N \widetilde{S_i}' (T_d' \operatorname{pr}^\perp(U) F_d) \widetilde{S_i} \right) \left( \sum_{i=1}^N \widetilde{S_i}' (F_d' \operatorname{pr}^\perp(U) F_d) \widetilde{S_i} \right)^- \\ &\quad \times \left( \sum_{i=1}^N \widetilde{S_i}' (F_d' \operatorname{pr}^\perp(U) T_d) \widetilde{S_i} \right). \end{aligned} \tag{3.4}$$

The equality sign in (3.4) will hold when (i) $l_{dik} = r_{di}/p$, $i = 0, \dots, t$, and (ii) $T_d' \operatorname{pr}^\perp(U) T_d$, $T_d' \operatorname{pr}^\perp(U) F_d$ and $F_d' \operatorname{pr}^\perp(U) F_d$ are invariant under any permutation of test treatments. Then by utilizing the definition of $\widetilde{S_i}$, we observe that $\sum_{i=1}^N \widetilde{S_i}'(T_d' \operatorname{pr}^\perp(U)T_d)\widetilde{S_i}$, $\sum_{i=1}^N \widetilde{S_i}'(T_d' \operatorname{pr}^\perp(U)F_d)\widetilde{S_i}$ and $\sum_{i=1}^N \widetilde{S_i}'(F_d' \operatorname{pr}^\perp(U)\times F_d)\widetilde{S_i}$ have the form

$$\begin{pmatrix} a_m & f_m J_{1\times t} \\ c_m J_{t\times 1} & (b_m - e_m)I_{t\times t} + e_m J_{t\times t} \end{pmatrix}, \qquad m = 1,2,3, \tag{3.5}$$

respectively. It can be shown that for $\sum_{i=1}^N \widetilde{S_i}'(T_d' \operatorname{pr}^\perp(U)T_d)\widetilde{S_i}$,

$$\begin{aligned}
a_1 &= N\left(r_{d0} - \frac{1}{p}\sum_{u=1}^n n_{d0u}^2\right),\\
b_1 &= \frac{N}{t}\left(np - r_{d0} - \frac{1}{p}\sum_{u=1}^n\sum_{i=1}^t n_{diu}^2\right),\\
c_1 &= f_1 = -\frac{a_1}{t},\\
e_1 &= -\frac{b_1}{t-1} + \frac{a_1}{t(t-1)};
\end{aligned}$$

for $\sum_{i=1}^N \widetilde{S_i}'(T_d' \operatorname{pr}^\perp(U)F_d)\widetilde{S_i}$,

$$\begin{aligned}
a_2 &= N\left(m_{d00} - \frac{1}{p}\sum_{u=1}^n n_{d0u}\tilde{n}_{d0u}\right),\\
b_2 &= \frac{N}{t}\sum_{i=1}^t\left(m_{dii} - \frac{1}{p}\sum_{u=1}^n n_{diu}\tilde{n}_{diu}\right),\\
c_2 &= -\frac{a_2}{t},\\
f_2 &= \frac{N}{t}\left(\hat{r}_{d0} - \frac{p-1}{p}r_{d0} - m_{d00} + \frac{1}{p}\sum_{u=1}^n n_{d0u}\tilde{n}_{d0u}\right),\\
e_2 &= -\frac{b_2+f_2}{t-1};
\end{aligned}$$

and for $\sum_{i=1}^N \widetilde{S_i}'(F_d' \operatorname{pr}^\perp(U)F_d)\widetilde{S_i}$,

$$\begin{aligned}
a_3 &= N\left(\tilde{r}_{d0} - \frac{1}{p}\sum_{u=1}^n \tilde{n}_{d0u}^2\right),\\
b_3 &= \frac{N}{t}\left(n(p-1) - \tilde{r}_{d0} - \frac{1}{p}\sum_{u=1}^n\sum_{i=1}^t \tilde{n}_{diu}^2\right),\\
c_3 &= f_3 = \frac{N}{t}\left(-\frac{p-1}{p}\tilde{r}_{d0} + \frac{1}{p}\sum_{u=1}^n \tilde{n}_{d0u}^2\right),
\end{aligned}$$

$$e_3 = \frac{N}{t(t-1)}\left\{\frac{p-1}{p}\tilde{r}_{d0} - \frac{1}{p}\sum_{u=1}^{n}\tilde{n}_{d0u}^2 \right.$$
$$\left. - \frac{p-1}{p}(n(p-1)-\tilde{r}_{d0}) + \frac{1}{p}\sum_{u=1}^{n}\sum_{i=1}^{t}\tilde{n}_{diu}^2\right\}.$$

Thus, it is easy to see that when the control treatment is kept in the same positions, $\sum_{i=1}^N \widetilde{S_i}'(F_d' \operatorname{pr}^\perp(U)F_d)\widetilde{S_i}$ will be maximized under the Loewner ordering when each test treatment appears at most once in each of $n$ first $p-1$ periods, that is,

$$\sum_{i=1}^N \widetilde{S_i}'(F_d' \operatorname{pr}^\perp(U)F_d)\widetilde{S_i} \le \begin{pmatrix} a_3 & c_3J_{1\times t} \\ c_3J_{t\times 1} & (\tilde{b}_3-\tilde{e}_3)I_{t\times t}+\tilde{e}_3J_{t\times t}\end{pmatrix}, \tag{3.6}$$

with $\tilde{b}_3 = \frac{N[n(p-1)-\tilde{r}_{d0}](p-1)}{tp}$ and $\tilde{e}_3 = \frac{N}{t(t-1)}\{\frac{p-1}{p}\tilde{r}_{d0} - \frac{1}{p}\sum_{u=1}^n \tilde{n}_{d0u}^2 - \frac{p-2}{p}(n\times(p-1)-\tilde{r}_{d0})\}$. When $0<\tilde{r}_{d0}<n(p-1)$, it can be easily verified that $a_3>0$, $\tilde{b}_3>\tilde{e}_3$ and $a_3\tilde{b}_3+(t-1)a_3\tilde{e}_3-tc_3^2>0$. So the right-hand side of (3.6) is a positive definite matrix and its inverse has the same matrix form as in (3.6) when $a_3$, $\tilde{b}_3$, $c_3$ and $\tilde{e}_3$ are replaced by $a_4$, $b_4$, $c_4$ and $e_4$. Here

$$a_4 = \frac{\tilde{b}_3+(t-1)\tilde{e}_3}{a_3\tilde{b}_3+(t-1)a_3\tilde{e}_3-tc_3^2},$$

$$b_4 = \frac{a_3\tilde{b}_3+(t-2)a_3\tilde{e}_3-(t-1)c_3^2}{(\tilde{b}_3-\tilde{e}_3)(a_3\tilde{b}_3+(t-1)a_3\tilde{e}_3-tc_3^2)},$$

$$c_4 = \frac{-c_3}{a_3\tilde{b}_3+(t-1)a_3\tilde{e}_3-tc_3^2},$$

$$e_4 = \frac{c_3^2-a_3\tilde{e}_3}{(\tilde{b}_3-\tilde{e}_3)(a_3\tilde{b}_3+(t-1)a_3\tilde{e}_3-tc_3^2)}.$$

The inverse matrix of the right-hand side of (3.6) can be expressed as $D+c_4J$, where

$$D = \begin{pmatrix} a_4-c_4 & 0_{1\times t} \\ 0_{t\times 1} & (b_4-e_4)I_{t\times t}+(e_4-c_4)J_{t\times t}\end{pmatrix}.$$

Notice that $c_4\ge 0$ due to $c_3\le 0$, and we have

$$\left(\sum_{i=1}^N \widetilde{S_i}'(T_d' \operatorname{pr}^\perp(U)F_d)\widetilde{S_i}\right)$$
$$\times\left(\sum_{i=1}^N \widetilde{S_i}'(F_d' \operatorname{pr}^\perp(U)F_d)\widetilde{S_i}\right)^-\left(\sum_{i=1}^N \widetilde{S_i}'(F_d' \operatorname{pr}^\perp(U)T_d)\widetilde{S_i}\right) \tag{3.7}$$

$$\geq \left(\sum_{i=1}^{N} \widetilde{S_i}'(T_d' \operatorname{pr}^{\perp}(U)F_d)\widetilde{S_i}\right) D \left(\sum_{i=1}^{N} \widetilde{S_i}'(F_d' \operatorname{pr}^{\perp}(U)T_d)\widetilde{S_i}\right).$$

The equality sign in the above inequality will hold when $l_{d01} = r_{d0}/p$ and each test treatment appears at most once in each of $n$ first $p-1$ periods. By (3.3), (3.4), (3.7) and by direct calculations, we have

$$M_d \leq \bar{M}_d \leq \widetilde{M}_d,$$

and the equality signs hold when the three conditions in Lemma 4 are satisfied. Here $\widetilde{M}_d = xI + yJ$, where

$$x = \frac{1}{N}[b_1 - e_1 - (b_2 - e_2)^2(b_4 - e_4)]$$

and

$$y = \frac{1}{N}[e_1 - c_2^2(a_4 - c_4) - e_2(b_2 - e_2 - f_2)(b_4 - e_4) - f_2^2(e_4 - c_4)].$$

Thus, to reach our conclusion, it is enough to show that the eigenvalues of $\widetilde{M}_d$ are $x_0$ with multiplicity $t-1$ and $y_0$ with multiplicity 1.

Since $\widetilde{M}_d = xI + yJ$, its eigenvalues are $x$ with multiplicity $t-1$ and $x + ty$ with multiplicity 1. It is easy to check that $b_4 - e_4 = 1/(\tilde{b}_3 - \tilde{e}_3)$, $b_2 - e_2 = (tb_2 + f_2)/(t-1)$ and $b_1 - e_1 = (t^2 b_1 - a_1)/(t(t-1))$. By direct calculations, we obtain $x = x_0$.

Now it suffices to show that $x + ty = y_0$. Notice that $b_2 + f_2 + (t-1)e_2 = 0$ and by the definitions of $a_4$, $b_4$, $c_4$ and $e_4$, we have

$$\begin{aligned} x + ty &= \frac{1}{N}[b_1 + (t-1)e_1 - tc_2^2(a_4 - c_4) - tf_2^2(e_4 - c_4) - f_2^2(b_4 - e_4)] \\ &= \frac{a_1}{Nt} - \frac{tc_2^2(\tilde{b}_3 + (t-1)\tilde{e}_3 + c_3) + f_2^2(a_3 + tc_3)}{N[a_3\tilde{b}_3 + (t-1)a_3\tilde{e}_3 - tc_3^2]}. \end{aligned}$$

By direct calculations, we obtain

$$\begin{aligned} a_3\tilde{b}_3 + (t-1)a_3\tilde{e}_3 - tc_3^2 &= \frac{N^2}{tp^2}\left[np(p-1)\tilde{r}_{d0} - \tilde{r}_{d0}^2 - n(p-1)\sum_{u=1}^{n} \tilde{n}_{d0u}^2\right], \\ \tilde{b}_3 + (t-1)\tilde{e}_3 + c_3 &= \frac{N}{tp}[n(p-1) - \tilde{r}_{d0}], \\ a_3 + tc_3 &= \frac{N\tilde{r}_{d0}}{p}. \end{aligned}$$

Thus $x + ty = y_0$. $\square$

**4. Optimal crossover designs for test treatments versus a control treatment.** In this section we will construct a family of crossover designs and prove their optimality for comparing $t$ test treatments versus a control treatment over a large class of designs $d$ in $\Omega_{t+1,n,p}$ for which $l_{d0k} = r_{d0}/p$, $k = 1, \dots, p$ (i.e., the control treatment appears equally often in $p$ periods) and $m_{dii} = 0$ for $0 \leq i \leq t$ (i.e., no treatment is allowed to be preceded by itself in $d$). The class of such designs is denoted by $\Lambda_{t+1,n,p}$. Before presenting the main theorem, we need some preliminary lemmas. To make our notation simple, we define the following four expressions:

$$\begin{aligned}\Delta_1 &= t(p-1)(np - r_{d0}) - p\left(r_{d0} - \frac{1}{p}\sum_{u=1}^{n} n_{d0u}^2\right)\\ &\quad - \frac{(nt(p-1) - t\tilde{r}_{d0} - \sum_{u=1}^{n} n_{d0u}\tilde{n}_{d0u})^2}{n(p-1)(pt-t-1) - (pt-t+p-2)\tilde{r}_{d0} + \sum_{u=1}^{n} \tilde{n}_{d0u}^2},\\ \Delta_2 &= p\left(r_{d0} - \frac{1}{p}\sum_{u=1}^{n} n_{d0u}^2\right) - \frac{n(p-1)(\sum_{u=1}^{n} n_{d0u}\tilde{n}_{d0u})^2}{np(p-1)\tilde{r}_{d0} - \tilde{r}_{d0}^2 - n(p-1)\sum_{u=1}^{n}\tilde{n}_{d0u}^2},\\ \Theta_1 &= \frac{nt(p-1) - t\tilde{r}_{d0} - \sum_{u=1}^{n} n_{d0u}\tilde{n}_{d0u}}{n(p-1)(pt-t-1) - (pt-t+p-2)\tilde{r}_{d0} + \sum_{u=1}^{n}\tilde{n}_{d0u}^2},\\ \Theta_2 &= \frac{n(p-1)\sum_{u=1}^{n} n_{d0u}\tilde{n}_{d0u}}{np(p-1)\tilde{r}_{d0} - \tilde{r}_{d0}^2 - n(p-1)\sum_{u=1}^{n}\tilde{n}_{d0u}^2}.\end{aligned}$$

The next lemma provides a simpler lower bound than that in Lemma 4 for any design $d$ in $\Lambda_{t+1,n,p}$. This lower bound depends only on the property of the control treatment in the design $d$.

LEMMA 5. *For any $d \in \Lambda_{t+1,n,p}$ with fixed $r_{d0}$, we have*

$$\operatorname{Tr}(M_d^{-1}) \geq \frac{t(t-1)^2 p}{x_1} + \frac{tp}{y_1}. \tag{4.1}$$

*Here*

$$\begin{aligned}x_1 &= t(p-1)(np - r_{d0}) - p\left(r_{d0} - \frac{1}{p}\min\sum_{u=1}^{n} n_{d0u}^2\right)\\ &\quad - \frac{[nt(p-1) - t\tilde{r}_{d0} - \min\sum_{u=1}^{n} n_{d0u}\tilde{n}_{d0u}]^2}{n(p-1)(pt-t-1) - (pt-t+p-2)\tilde{r}_{d0} + \min\sum_{u=1}^{n}\tilde{n}_{d0u}^2}\end{aligned}$$

*and*

$$y_1 = p\left(r_{d0} - \frac{1}{p}\min\sum_{u=1}^{n} n_{d0u}^2\right) - \frac{n(p-1)(\min\sum_{u=1}^{n} n_{d0u}\tilde{n}_{d0u})^2}{np(p-1)\tilde{r}_{d0} - \tilde{r}_{d0}^2 - n(p-1)\sum_{u=1}^{n}\tilde{n}_{d0u}^2}.$$

*The equality in (4.1) will hold when the following two conditions are satisfied in addition to the three conditions in Lemma 4:*

(i) $n_{diu}$, $i=1,\dots,t, u=1,\dots,n$, *are either* 0 *or* 1,
(ii) $\sum_{u=1}^n n_{d0u}^2$, $\sum_{u=1}^n n_{d0u}\tilde{n}_{d0u}$ *and* $\sum_{u=1}^n \tilde{n}_{d0u}^2$ *are minimized.*

*Here, the minimizations of* $\sum_{u=1}^n n_{d0u}^2$, $\sum_{u=1}^n n_{d0u}\tilde{n}_{d0u}$ *and* $\sum_{u=1}^n \tilde{n}_{d0u}^2$ *are over all possible designs* $d \in \Lambda_{t+1,n,p}$ *with fixed* $r_{d0}$.

PROOF. Since for any design $d \in \Lambda_{t+1,n,p}$, with $r_{d0}=0$ or $np$, the contrasts between test treatments and the control treatment are not estimable, thus only designs with $0 < \tilde{r}_{d0} < n(p-1)$ are considered here. When $0 < r_{d0} < np$, from Lemma 4, we know that $\mathrm{Tr}(M_d^{-1}) \geq \frac{t-1}{x_0} + \frac{1}{y_0}$, with $x_0$ and $y_0$ as defined in Lemma 4. By Proposition A.2 (for this and all other propositions see the Appendix), we notice that $t\sum_{i=1}^t \sum_{u=1}^n n_{diu}\tilde{n}_{diu} \geq t(n(p-1)-\tilde{r}_{d0}) \geq \sum_{u=1}^n n_{d0u}\tilde{n}_{d0u}$. We also know that $m_{dii}=0$, $i=0,\dots,t$, and $\hat{r}_{d0} = \frac{(p-1)r_{d0}}{p}$. Thus

$$x_0 \leq \frac{1}{t(t-1)p}\Delta_1,$$

$$y_0 = \frac{1}{tp}\Delta_2,$$

and equality holds when $n_{diu}$, $i=1,\dots,t$, $u=1,\dots,n$, are either 0 or 1. So we have

$$\mathrm{Tr}(M_d^{-1}) \geq \frac{t(t-1)^2 p}{\Delta_1} + \frac{tp}{\Delta_2}. \tag{4.2}$$

By Proposition A.2, we have

$$\sum_{u=1}^n n_{d0u}\tilde{n}_{d0u} \leq nt(p-1) - t\tilde{r}_{d0}.$$

For any $d \in \Lambda_{t+1,n,p}$ with fixed $r_{d0}$, let $\xi_1$, $\xi_2$ and $\xi_3$ be independent variables satisfying the following inequalities:

$$\xi_1 \geq \min_d\left(\sum_{u=1}^n n_{d0u}^2\right),$$

$$\min_d\left(\sum_{u=1}^n n_{d0u}\tilde{n}_{d0u}\right) \leq \xi_2 \leq nt(p-1) - t\tilde{r}_{d0},$$

$$\xi_3 \geq \min_d\left(\sum_{u=1}^n \tilde{n}_{d0u}^2\right).$$

By Proposition A.6, we know that $\min_d(\sum_{u=1}^n n_{d0u}^2)$, $\min_d(\sum_{u=1}^n n_{d0u}\tilde{n}_{d0u})$ and $\min_d(\sum_{u=1}^n \tilde{n}_{d0u}^2)$ are fixed for given $n$, $p$ and $r_{d0}$.

Define $\widetilde{\Delta}_1$, $\widetilde{\Delta}_2$, $\widetilde{\Theta}_1$ and $\widetilde{\Theta}_2$ to be the corresponding values of $\Delta_1$, $\Delta_2$, $\Theta_1$ and $\Theta_2$, respectively, after replacing $\sum_{u=1}^n n_{d0u}^2$, $\sum_{u=1}^n n_{d0u}\tilde{n}_{d0u}$ and $\sum_{u=1}^n \tilde{n}_{d0u}^2$ by $\xi_1$, $\xi_2$ and $\xi_3$, respectively, in $\Delta_1$, $\Delta_2$, $\Theta_1$ and $\Theta_2$. Note that $\sum_{u=1}^n n_{d0u}^2$, $\sum_{u=1}^n n_{d0u}\tilde{n}_{d0u}$ and $\sum_{u=1}^n \tilde{n}_{d0u}^2$ are not independent of each other, and thus the set $\{\sum_{u=1}^n n_{d0u}^2, \sum_{u=1}^n n_{d0u}\tilde{n}_{d0u}, \sum_{u=1}^n \tilde{n}_{d0u}^2\}$ is just a subset of $\{\xi_1, \xi_2, \xi_3\}$. Thus for any $d \in \Lambda_{t+1,n,p}$ with fixed $r_{d0}$, we have

$$\min_{\sum_{u=1}^n n_{d0u}^2, \sum_{u=1}^n n_{d0u}\tilde{n}_{d0u}, \sum_{u=1}^n \tilde{n}_{d0u}^2} \left(\frac{t(t-1)^2p}{\Delta_1} + \frac{tp}{\Delta_2}\right) \geq \min_{\xi_1,\xi_2,\xi_3} \left(\frac{t(t-1)^2p}{\widetilde{\Delta}_1} + \frac{tp}{\widetilde{\Delta}_2}\right).$$

Define

$$H(r_{d0}, \xi_1, \xi_2, \xi_3) = \frac{t(t-1)^2p}{\widetilde{\Delta}_1} + \frac{tp}{\widetilde{\Delta}_2}.$$

To reach our conclusion, it is sufficient to show that for fixed $r_{d0}$, $H(r_{d0}, \xi_1, \xi_2, \xi_3)$ is an increasing function of $\xi_1$, $\xi_2$ and $\xi_3$.

By direct calculation and Proposition A.3, we have

$$\frac{\partial H(r_{d0}, \xi_1, \xi_2, \xi_3)}{\partial \xi_1} = tp\left(\frac{1}{\widetilde{\Delta}_2^2} - \frac{(t-1)^2}{\widetilde{\Delta}_1^2}\right) \geq 0.$$

Also we have

$$\frac{\partial H(r_{d0}, \xi_1, \xi_2, \xi_3)}{\partial \xi_2} = 2tp\left(\frac{\widetilde{\Theta}_2}{\widetilde{\Delta}_2^2} - \frac{(t-1)^2\widetilde{\Theta}_1}{\widetilde{\Delta}_1^2}\right) \geq 0. \tag{4.3}$$

The last inequality holds when $\frac{\tilde{r}_{d0}}{n} > \frac{p-1}{t+1}$ due to Proposition A.3 and (A.10) in Proposition A.5. Further, by applying Proposition A.4 and (A.9) in Proposition A.5, the inequality still holds when $\frac{\tilde{r}_{d0}}{n} \leq \frac{p-1}{t+1}$.

By direct calculations, we have

$$\frac{\partial H(r_{d0}, \xi_1, \xi_2, \xi_3)}{\partial \xi_3} = tp\left(\frac{\widetilde{\Theta}_2^2}{\widetilde{\Delta}_2^2} - \frac{(t-1)^2\widetilde{\Theta}_1^2}{\widetilde{\Delta}_1^2}\right). \tag{4.4}$$

And by the same argument, we have $\frac{\partial H(r_{d0},\xi_1,\xi_2,\xi_3)}{\partial \xi_3} \geq 0$. Therefore for a given $r_{d0}$, $\frac{t(t-1)^2p}{\widetilde{\Delta}_1} + \frac{tp}{\widetilde{\Delta}_2}$ will achieve the minimum value when $\xi_1$, $\xi_2$ and $\xi_3$ are minimized. Thus, the conclusion is obtained. □

Next we will introduce some definitions which are similar to the definitions in Kunert and Stufken (2002). A design $d \in \Omega_{t+1,n,p}$ is called:

(i) a balanced test-control incomplete block design for the direct effects (with units as blocks) if (a) each test treatment $1 \leq i \leq t$ appears equally often in the design; (b) each test treatment appears in each unit at most once; (c) the number of units where any two test treatments $i$ and $j$ both appear is the same for every $i \neq j$, $1 \leq i \leq t$, $1 \leq j \leq t$; (d) the control treatment appears in each unit either $[\frac{r_{d0}}{n}]$ or $[\frac{r_{d0}}{n}]+1$ times; and (e) the number of units where the control treatment appears $[\frac{r_{d0}}{n}]$ times and test treatment $i$ appears is the same for every $1 \le i \le t$.

(ii) a balanced test-control incomplete block design for the carryover effects (with units as blocks) if the first $p-1$ periods of the design form a balanced test-control block design for the direct effects in $\Omega_{t+1,n,p-1}$.

(iii) a balanced for test-control carryover effects design if (a) every test treatment $1 \le i \le t$ is immediately preceded by every other test treatment $1 \le j \le t$ equally often for every $i \ne j$; (b) the control treatment is immediately preceded by every test treatment $1 \le i \le t$ equally often; (c) every test treatment $1 \le i \le t$ is immediately preceded by the control treatment equally often, and no treatment including the control is immediately preceded by itself.

(iv) a proportional frequency design for test-control on the periods if every test treatment $1 \le i \le t$ appears in every period exactly $\frac{np-r_{d0}}{tp}$ times, and the control treatment appears in every period exactly $\frac{r_{d0}}{p}$ times.

A design $d^* \in \Omega_{t+1,n,p}$ is called a totally balanced test-control incomplete crossover design if:

(i) $d^*$ is a balanced test-control incomplete block design for the direct effects,

(ii) $d^*$ is a balanced test-control incomplete block design for the carryover effects,

(iii) $d^*$ is balanced for test-control carryover effects,

(iv) $d^*$ is a proportional frequency design for test-control on the periods,

(v) the number of units where test treatment $j$ appears once in the first $p-1$ periods and test treatment $i$ appears once is the same for every pair $i \ne j$; the number of units where the control treatment appears $[\frac{\tilde{r}_{d0}}{n}]$ times in the first $p-1$ periods and test treatment $i$ appears once is the same for every test treatment $1 \le i \le t$; the number of units where test treatment $i$ appears once in the first $p-1$ periods and the control treatment appears $[\frac{r_{d0}}{n}]$ times is the same for every test treatment $1 \le i \le t$.

Examples of such designs will be given in Section 5. The following lemma summarizes some useful properties of a totally balanced test-control incomplete crossover design.

LEMMA 6. *If $d$ is a totally balanced test-control incomplete crossover design, then:*

(i) $l_{dik} = r_{di}/p$, $i = 0, \ldots, t$,

(ii) *$M_d$ is a completely symmetric matrix [since $T_d' \operatorname{pr}^{\perp}(U) T_d$, $T_d' \operatorname{pr}^{\perp}(U) F_d$ and $F_d' \operatorname{pr}^{\perp}(U) F_d$ are invariant under any permutation of test treatments],*

(iii) $\sum_{u=1}^{n} n_{d0u}^2$, $\sum_{u=1}^{n} n_{d0u}\tilde{n}_{d0u}$ *and* $\sum_{u=1}^{n} \tilde{n}_{d0u}^2$ *are minimized over the designs with fixed* $r_{d0}$.

We are now ready to summarize our findings in the following theorem.

THEOREM 1. *For* $p \le t+1$, *a design* $d^*$ *is simultaneously A-optimal and MV-optimal in* $\Lambda_{t+1,n,p}$ *if* $d^*$ *is a totally balanced test-control incomplete crossover design and* $r_{d^*0}$ *minimize the right-hand side of* (4.1).

PROOF. By Lemma 6, the conditions for the equality sign in Lemma 5 hold. From Proposition A.6, we notice that $\min \sum_{u=1}^{n} n_{d0u}^2$, $\min \sum_{u=1}^{n} n_{d0u}\tilde{n}_{d0u}$ and $\min \sum_{u=1}^{n} \tilde{n}_{d0u}^2$ are functions of $r_{d0}$. Thus we have

$$\mathrm{Tr}(M_{d^*}^{-1}) = \frac{t(t-1)^2 p}{x_1^*} + \frac{tp}{y_1^*} = \min_{\tilde{r}_{d0}} \left( \frac{t(t-1)^2 p}{x_1} + \frac{tp}{y_1} \right) \le \mathrm{Tr}(M_d^{-1}),$$

where $x_1^*$ and $y_1^*$ are the corresponding $x_1$ and $y_1$ when $d$ is $d^*$. The last inequality holds due to Lemma 5. Thus design $d^*$ is A-optimal. And since $M_{d^*}^{-1}$ is a completely symmetric matrix by Lemma 6, the MV-optimality of design $d^*$ follows by Lemma 1. □

**5. Assorted mathematical tools useful for the construction of optimal crossover designs.** While Theorem 1 specifies a set of sufficient conditions for a crossover design to be simultaneously A-optimal and MV-optimal, it will not provide any mathematical tools for constructing these designs. The purpose of this section is to alleviate this deficiency and present some mathematical tools which in conjunction with some simple computer programs can help in constructing these designs. In practice, $t$ and $p$ are given and the job of the statistician is to find an optimal design for a given desirable $n$. If we want to rely on Theorem 1, then we are left with two tasks: first, to apply (A.13)–(A.15) in Proposition A.6 into the right-hand side of (4.1) and run a simple computer program to find the $r_{d^*0}$ which minimizes the right-hand side of (4.1), and second, to build a totally balanced test-control incomplete crossover design $d$ where the replication of the control treatment in this design is $r_{d^*0}$.

Our experience shows that the family of totally balanced test-control incomplete crossover designs based on $t$, $p$ and $n$ contains very useful designs for our problems if $r_{d0} = n$. We have observed that these designs are either optimal or at least highly efficient. Thus we shall concentrate our effort on this class of crossover designs. Not to leave the impression that the optimal designs are exclusively in this class of designs, we shall at the end of this section exhibit two optimal designs based on Theorem 1 with $r_{d0} \ne n$. But, we should point out that the construction of optimal designs with $r_{d0} \ne n$ is not easy and we need more mathematical tools in this area.

LEMMA 7. *For given $t$, $p$, $n$ and $r_{d0} = n$, a totally balanced test-control incomplete crossover design exists only if both $\frac{n(p-1)}{pt}$ and $\frac{(p-1)(p-2)n}{pt(t-1)}$ are integers.*

PROOF. Notice that the condition $\frac{n(p-1)}{pt}$ is an integer implies $\frac{n}{p}$ is also an integer. This condition is necessary due to the uniformity for the test treatments on the periods. The condition that $\frac{(p-1)(p-2)n}{pt(t-1)}$ is an integer is also necessary since the design is balanced for test treatments. In fact, if the control treatment appears in the first or the last period, then any test treatment will be preceded by every other test treatment $p-2$ times. And if the control treatment appears in any other periods, then any test treatment will be preceded by every other test treatment $p-3$ times. So, totally we have $(p-1)(p-2)\frac{n}{p}$ times that one test treatment is preceded by every other test treatment. Thus $\frac{(p-1)(p-2)n}{pt(t-1)}$ must be an integer. □

As we shall see soon, for many cases the necessary conditions in Lemma 7 become sufficient for the existence of totally balanced test-control incomplete crossover designs. Lemmas 8 and 9 deal with cases for which $p = t + 1$ and $p = t$. As for the case of $p < t$, we shall provide a very useful tool that could help in developing such crossover designs.

LEMMA 8. *For given $t$, $p = t + 1$, $n$ and $r_{d0} = n$, a totally balanced test-control incomplete crossover design exists if there is a balanced uniform design in $\Omega_{t+1,n,p}$.*

PROOF. Simply replace treatment $t + 1$ in the balanced uniform design by the control treatment. Then it is easy to see that the modified design is a totally balanced test-control incomplete crossover design with the stated parameters. □

Notice that a necessary condition for the existence of a balanced uniform design in $\Omega_{t+1,n,t+1}$ is $n = \lambda(t + 1)$ for some positive integer $\lambda$. According to Higham (1998), the class $\Omega_{t+1,n,t+1}$ contains a balanced uniform design when either $n$ is an even multiple of $t + 1$ or $t + 1$ is a composite number, that is, it can be written as a product of two positive integers each larger than 1.

LEMMA 9. *For given $t$, $p = t$, $n$ and $r_{d0} = n$, a totally balanced test-control incomplete crossover design exists if* (i) *$\frac{n}{t^2}$ is an integer and $t$ is a composite number, or* (ii) *$\frac{n}{t^2}$ is an even integer and $t$ is a prime number.*

PROOF. We shall give an explicit way of constructing these designs. By the two conditions in the lemma, $\frac{n}{t}$ is a multiple of $t$ when $t$ is a composite number and $\frac{n}{t}$ is an even multiple of $t$ when $t$ is a prime number. Thus we can always construct a balanced uniform design in $\Omega_{t,n/t,t}$. The totally balanced test-control incomplete crossover design with the stated parameters can be constructed in the following way: (i) Construct a balanced uniform design in $\Omega_{t,n/t,t}$. (ii) Replace treatment 1 by the control treatment 0 in the balanced uniform design. (iii) Repeat step (ii) for the remaining test treatments to produce in total $t$ new designs. (iv) Juxtapose the $t$ new designs into a single design. The resulting design is the desired design. □

For the case $p < t$, our experience shows that the following steps lead successfully to a totally balanced test-control incomplete crossover design.

STEP 1. Construct a balanced incomplete block (BIB) design based on $t$ test treatments in blocks of size $p-1$. Let $b$ be the number of blocks in this BIB design and number the blocks from 1 to $b$ in an arbitrary fashion.

STEP 2. Construct $p$ arrays each of size $p$ by $b$. Fill the $k$th array, $k = 1, 2, \ldots, p$, in the following way. Fill the entire $b$ cells in the $k$th row of this array by the control treatment. Note that each of the $b$ columns of the array now has $p-1$ empty cells. Fill these $b$ columns arbitrarily with the $b$ blocks of the BIB design. In this way we have produced $p$ arrays each filled with the control treatment or the blocks of the BIB design. Juxtapose these $p$ arrays into a big $p$ by $pb$ array.

STEP 3. Shuffle the positions of the test treatments in each column of the array produced in Step 2, for the purpose of converting the $p$ by $pb$ array into a totally balanced test-control incomplete crossover design.

Although the preceding three steps cannot guarantee that in all cases a totally balanced test-control incomplete crossover design will be produced, our experience indicates that they are a very useful mathematical tool in building such designs when they exist. Indeed, we have succeeded in constructing all small totally balanced test-control incomplete crossover designs or designs which are very close to such crossover designs.

In the sequel we present some interesting examples based on the tools presented here. But first, a brief overview. For given $t$ and $p \leq t+1$ we can choose $n$ so that the two integer conditions in Lemma 7 are satisfied. Then we can apply (A.13)–(A.15) in Proposition A.6 into the right-hand side of (4.1) and run a simple computer program to find the $r_{d^*0}$ which minimizes the right-hand side of (4.1). If $r_{d^*0} = n$, we can utilize Lemmas 8 and 9 as the steps for the cases $p < t$ and construct a totally balanced test-control incomplete crossover design which is simultaneously A-optimal and MV-optimal by Theorem 1. In case $r_{d^*0} \neq n$, we can compare the corresponding minimum value of the right-hand side of (4.1) with the value when $r_{d0} = n$. If the two values are very close, we can use Lemmas 8 and 9 as the steps for $p < t$ to build a highly efficient crossover design. Ideally, we should investigate to see if it is possible to construct an optimal design when $r_{d^*0} \neq n$.

We shall now concentrate on finding small size optimal designs for these practically desirable cases, namely (i) $p = 3$, (ii) $p = 4$ and (iii) $p = 5$. In addition, we exhibit two optimal crossover designs when $r_{d^*0} \neq n$. The general method of constructing optimal crossover designs when $r_{d^*0} \neq n$ is very difficult and remains open.

5.1. *Simultaneous A-optimal and MV-optimal crossover designs for three periods.* When $t = 2$, 3 or 4, then $n$ must be a multiple of 3, 9 or 18, respectively. For the minimum values of $n$ which satisfy the integer conditions in Lemma 7, we found that $r_{d^*0} = n$. Unfortunately, there are no totally balanced test-control incomplete crossover designs for either $t = 2$ and $n = 3$ or $t = 4$ and $n = 18$. Thus, we searched for optimal designs for the next allowable value of $n$. We tried $n = 6$ for $t = 2$ and $n = 36$ for $t = 4$. Fortunately, $r_{d^*0} = 6$ and 36, respectively, for these situations and both designs can be easily constructed. For $t = 3$, the minimum $n$ which satisfies the two integer conditions in Lemma 7 is 9. And indeed for this case $r_{d^*0} = 9$. Example 1 exhibits one such optimal design for these parameters.

Example 1. A-optimal and MV-optimal design for $p = 3$, $t = 3$ and $n = 9$:

| | | | | | | | | |
|---|---|---|---|---|---|---|---|---|
| 0 | 0 | 0 | 2 | 3 | 1 | 2 | 3 | 1 |
| 1 | 2 | 3 | 0 | 0 | 0 | 1 | 2 | 3 |
| 2 | 3 | 1 | 1 | 2 | 3 | 0 | 0 | 0 |

When $t = 5$, $n$ must be a multiple of 30 and indeed $n = 30$ satisfies the integer conditions in Lemma 7 and $r_{d^*0} = 30$. Example 2 exhibits one such optimal design for these parameters.

Example 2. A-optimal and MV-optimal design for $p = 3$, $t = 5$ and $n = 30$:

| | | | | | | | | | | | | | | | | | | | | | | | | | | | | | |
|---|---|---|---|---|---|---|---|---|---|---|---|---|---|---|---|---|---|---|---|---|---|---|---|---|---|---|---|---|---|
| 0 | 0 | 0 | 0 | 0 | 0 | 0 | 0 | 0 | 0 | 1 | 1 | 4 | 5 | 3 | 2 | 2 | 4 | 3 | 5 | 1 | 1 | 4 | 5 | 3 | 2 | 2 | 4 | 3 | 5 |
| 2 | 3 | 1 | 1 | 2 | 4 | 5 | 3 | 5 | 4 | 0 | 0 | 0 | 0 | 0 | 0 | 0 | 0 | 0 | 0 | 2 | 3 | 1 | 1 | 2 | 4 | 5 | 3 | 5 | 4 |
| 1 | 1 | 4 | 5 | 3 | 2 | 2 | 4 | 3 | 5 | 2 | 3 | 1 | 1 | 2 | 4 | 5 | 3 | 5 | 4 | 0 | 0 | 0 | 0 | 0 | 0 | 0 | 0 | 0 | 0 |

5.2. *Simultaneous A-optimal and MV-optimal crossover designs for four periods.* We tried to construct these designs for $t = 3$, 4, 5, 6, 7, 8 and 9. The story is as follows. First, we note that the minimum $n$ for these values of $t$ must be a multiple of 4, 16, 40, 40, 28, 224 and 48, respectively. So we tried the minimum $n$. Fortunately, for all these cases $r_{d^*0} = n$. We did not attempt to construct the design for the case of $t = 8$ and $n = 224$ due to the size of $n$. For the remaining six cases we succeeded in constructing optimal crossover designs in the form of Theorem 1. For $t = 3$ and $n = 4$, we can apply the technique in Lemma 8. When $t = 4$ and $n = 16$, we can use the technique of Lemma 9. And for $t = 5$ and $n = 40$, $t = 6$ and $n = 40$, $t = 7$ and $n = 28$ and $t = 9$ and $n = 48$, we successfully used the construction steps given after Lemma 9. We give samples of such designs in Examples 3–6.

EXAMPLE 3. A-optimal and MV-optimal design for $p = 4$, $t = 5$ and $n = 40$:

0 0 0 0 0 0 0 0 0 0 1 3 5 2 4 2 5 3 1 4 1 3 5 2 4 2 5 3 1 4 1 2 3 4 5 4 3 2 1 5
1 2 3 4 5 4 3 2 1 5 0 0 0 0 0 0 0 0 0 0 4 1 3 5 2 4 2 5 3 1 5 1 2 3 4 5 4 3 2 1
5 1 2 3 4 5 4 3 2 1 4 1 3 5 2 4 2 5 3 1 0 0 0 0 0 0 0 0 0 0 2 3 4 5 1 3 2 1 5 4
2 3 4 5 1 3 2 1 5 4 3 5 2 4 1 5 3 1 4 2 3 5 2 4 1 5 3 1 4 2 0 0 0 0 0 0 0 0 0 0

EXAMPLE 4. A-optimal and MV-optimal design for $p = 4$, $t = 6$ and $n = 40$:

0 0 0 0 0 0 0 0 0 0 1 4 5 4 6 2 3 5 3 1 3 5 3 2 1 6 1 2 6 4 4 2 6 6 5 1 2 3 4 5
3 5 3 2 1 6 1 2 6 4 0 0 0 0 0 0 0 0 0 0 4 2 6 6 5 1 2 3 4 5 1 4 5 4 6 2 3 5 3 1
1 4 5 4 6 2 3 5 3 1 4 2 6 6 5 1 2 3 4 5 0 0 0 0 0 0 0 0 0 0 3 5 3 2 1 6 1 2 6 4
4 2 6 6 5 1 2 3 4 5 3 5 3 2 1 6 1 2 6 4 1 4 5 4 6 2 3 5 3 1 0 0 0 0 0 0 0 0 0 0

EXAMPLE 5. A-optimal and MV-optimal design for $p = 4$, $t = 7$ and $n = 28$:

0 0 0 0 0 0 0 3 4 5 6 7 1 2 1 2 3 4 5 6 7 4 5 6 7 1 2 3
1 2 3 4 5 6 7 0 0 0 0 0 0 0 4 5 6 7 1 2 3 3 4 5 6 7 1 2
3 4 5 6 7 1 2 4 5 6 7 1 2 3 0 0 0 0 0 0 0 1 2 3 4 5 6 7
4 5 6 7 1 2 3 1 2 3 4 5 6 7 3 4 5 6 7 1 2 0 0 0 0 0 0 0

EXAMPLE 6. A-optimal and MV-optimal design for $p = 4$, $t = 9$ and $n = 48$:

0 0 0 0 0 0 0 0 0 0 0 0 6 9 3 1 7 4 2 5 8 4 5 6 8 2 5 9 6 3 1 4 7 1 2 3 1 4 7 5 2 8 3 6 9 7 8 9
8 2 5 9 6 3 1 4 7 1 2 3 0 0 0 0 0 0 0 0 0 0 0 0 1 4 7 5 2 8 3 6 9 7 8 9 6 9 3 1 7 4 2 5 8 4 5 6
6 9 3 1 7 4 2 5 8 4 5 6 1 4 7 5 2 8 3 6 9 7 8 9 0 0 0 0 0 0 0 0 0 0 0 0 8 2 5 9 6 3 1 4 7 1 2 3
1 4 7 5 2 8 3 6 9 7 8 9 8 2 5 9 6 3 1 4 7 1 2 3 6 9 3 1 7 4 2 5 8 4 5 6 0 0 0 0 0 0 0 0 0 0 0 0

5.3. *Simultaneous A-optimal and MV-optimal crossover designs for five periods.* We tried to build such designs for $t=4$, 5, 6 and 7. For $t=4$, $n$ must be a multiple of 5 and the first $n$ for which a totally balanced test-control incomplete crossover design with $r_{d0}=n$ exists is $n=10$. For this value of $n$, we found that $r_{d^*0}=n$. The corresponding optimal design can be easily constructed using the tools in Lemma 8. However, we like to point out that if we need a design with bigger $n$ for this case, then $r_{d^*0}$ may not be $n$. In Example 10, we have exhibited an optimal design with $n=48$ and $r_{d^*0}=60$.

When $t=5$, $n$ must be a multiple of 25. So far we have not been able to construct a totally balanced test-control incomplete crossover design with $r_{d0}=n$ when $n=25$. However, we can construct such a design when $n=50$ by using the tools in Lemma 8. Notice that when $n=50$, $r_{d^*0}=60$ minimizes the right-hand side of (4.1) and the corresponding minimum value is 0.24179. Although $r_{d0}=50$ does not minimize the right-hand side of (4.1), its corresponding value is 0.24299, which is 99.5% efficient relative to the minimum value. So this design is highly efficient or even optimal.

When $t=6$, a totally balanced test-control incomplete crossover design will be relatively very large since $n$ must be multiple of 75. We tried $n=30$ knowing that we cannot use Theorem 1 to conclude optimality, but we hoped for a very efficient design. For $n=30$, we found $r_{d^*0}=30$ and the corresponding minimum value of the right-hand side of (4.1) is 0.55044. We used the construction steps after Lemma 9 and found a highly efficient design. This design is given in Example 7. For this design $d$, its $\mathrm{Tr}(M_d^{-1})=0.55419$, which is 99.3% efficient relative to the minimum value. So this design is highly efficient or even optimal.

EXAMPLE 7. Efficient design for $p=5$, $t=6$ and $n=30$ with efficiency $=99.3\%$:

| | | | | | | | | | | | | | | | | | | | | | | | | | | | | | |
|---|---|---|---|---|---|---|---|---|---|---|---|---|---|---|---|---|---|---|---|---|---|---|---|---|---|---|---|---|---|
| 0 | 0 | 0 | 0 | 0 | 0 | 6 | 1 | 2 | 3 | 4 | 5 | 2 | 3 | 4 | 5 | 6 | 1 | 3 | 4 | 5 | 6 | 1 | 2 | 5 | 6 | 1 | 2 | 3 | 4 |
| 1 | 2 | 3 | 4 | 5 | 6 | 0 | 0 | 0 | 0 | 0 | 0 | 6 | 1 | 2 | 3 | 4 | 5 | 5 | 6 | 1 | 2 | 3 | 4 | 6 | 1 | 2 | 3 | 4 | 5 |
| 6 | 1 | 2 | 3 | 4 | 5 | 2 | 3 | 4 | 5 | 6 | 1 | 0 | 0 | 0 | 0 | 0 | 0 | 2 | 3 | 4 | 5 | 6 | 1 | 2 | 3 | 4 | 5 | 6 | 1 |
| 2 | 3 | 4 | 5 | 6 | 1 | 5 | 6 | 1 | 2 | 3 | 4 | 3 | 4 | 5 | 6 | 1 | 2 | 0 | 0 | 0 | 0 | 0 | 0 | 1 | 2 | 3 | 4 | 5 | 6 |
| 5 | 6 | 1 | 2 | 3 | 4 | 3 | 4 | 5 | 6 | 1 | 2 | 4 | 5 | 6 | 1 | 2 | 3 | 6 | 1 | 2 | 3 | 4 | 5 | 0 | 0 | 0 | 0 | 0 | 0 |

When $t=7$, $n$ must be a multiple of 35. So far we have not succeeded in constructing an optimal design in the form of a totally balanced test-control incomplete crossover design with $r_{d0}=n=35$. However, we have been able to construct such a design for $n=70$. For $n=70$, we found $r_{d^*0}$ to be 70 and we shall present one such design in Example 8. The design in Example 8 is split into two parts $d_1$ and $d_2$ each based on $n=35$. While the union of the two designs is optimal, each of them is a highly efficient design in the class of designs with $n = 35$. When $n = 35$, the minimum value of the right-hand side of (4.1) is 0.61529 while $\mathrm{Tr}(M_{d_1}^{-1}) = 0.61904$ and $\mathrm{Tr}(M_{d_2}^{-1}) = 0.61927$. Thus $d_1$ is 99.39% efficient and $d_2$ is 99.36% efficient.

EXAMPLE 8. A-optimal and MV-optimal design for $p = 5$, $t = 7$ and $n = 70$:

$$d = d_1 \cup d_2,$$

where

$$d_1 = \begin{array}{c}
0\ 0\ 0\ 0\ 0\ 0\ 0\ 6\ 7\ 1\ 2\ 3\ 4\ 5\ 2\ 3\ 4\ 5\ 6\ 7\ 1\ 7\ 1\ 2\ 3\ 4\ 5\ 6\ 6\ 7\ 1\ 2\ 3\ 4\ 5 \\
2\ 3\ 4\ 5\ 6\ 7\ 1\ 0\ 0\ 0\ 0\ 0\ 0\ 0\ 7\ 1\ 2\ 3\ 4\ 5\ 6\ 6\ 7\ 1\ 2\ 3\ 4\ 5\ 7\ 1\ 2\ 3\ 4\ 5\ 6 \\
6\ 7\ 1\ 2\ 3\ 4\ 5\ 7\ 1\ 2\ 3\ 4\ 5\ 6\ 0\ 0\ 0\ 0\ 0\ 0\ 0\ 5\ 6\ 7\ 1\ 2\ 3\ 4\ 5\ 6\ 7\ 1\ 2\ 3\ 4 \\
5\ 6\ 7\ 1\ 2\ 3\ 4\ 2\ 3\ 4\ 5\ 6\ 7\ 1\ 5\ 6\ 7\ 1\ 2\ 3\ 4\ 0\ 0\ 0\ 0\ 0\ 0\ 0\ 2\ 3\ 4\ 5\ 6\ 7\ 1 \\
7\ 1\ 2\ 3\ 4\ 5\ 6\ 5\ 6\ 7\ 1\ 2\ 3\ 4\ 6\ 7\ 1\ 2\ 3\ 4\ 5\ 2\ 3\ 4\ 5\ 6\ 7\ 1\ 0\ 0\ 0\ 0\ 0\ 0\ 0
\end{array}$$

and

$$d_2 = \begin{array}{c}
0\ 0\ 0\ 0\ 0\ 0\ 0\ 6\ 7\ 1\ 2\ 3\ 4\ 5\ 5\ 6\ 7\ 1\ 2\ 3\ 4\ 7\ 1\ 2\ 3\ 4\ 5\ 6\ 6\ 7\ 1\ 2\ 3\ 4\ 5 \\
5\ 6\ 7\ 1\ 2\ 3\ 4\ 0\ 0\ 0\ 0\ 0\ 0\ 0\ 7\ 1\ 2\ 3\ 4\ 5\ 6\ 6\ 7\ 1\ 2\ 3\ 4\ 5\ 7\ 1\ 2\ 3\ 4\ 5\ 6 \\
6\ 7\ 1\ 2\ 3\ 4\ 5\ 7\ 1\ 2\ 3\ 4\ 5\ 6\ 0\ 0\ 0\ 0\ 0\ 0\ 0\ 2\ 3\ 4\ 5\ 6\ 7\ 1\ 2\ 3\ 4\ 5\ 6\ 7\ 1 \\
2\ 3\ 4\ 5\ 6\ 7\ 1\ 5\ 6\ 7\ 1\ 2\ 3\ 4\ 2\ 3\ 4\ 5\ 6\ 7\ 1\ 0\ 0\ 0\ 0\ 0\ 0\ 0\ 5\ 6\ 7\ 1\ 2\ 3\ 4 \\
7\ 1\ 2\ 3\ 4\ 5\ 6\ 2\ 3\ 4\ 5\ 6\ 7\ 1\ 6\ 7\ 1\ 2\ 3\ 4\ 5\ 5\ 6\ 7\ 1\ 2\ 3\ 4\ 0\ 0\ 0\ 0\ 0\ 0\ 0
\end{array}$$

Before closing this section we shall present two totally balanced test-control incomplete crossover designs for which $r_{d^*0} \neq n$. Such designs are extremely difficult to construct. For $p = 3$, $t = 7$ and $n = 49$, we have $r_{d^*0} = 42$ and for $p = 5$, $t = 4$ and $n = 48$, we have $r_{d^*0} = 60$.

EXAMPLE 9. A-optimal and MV-optimal design for $p = 3$, $t = 7$ and $n = 49$:

$$\begin{array}{c}
0\ 0\ 0\ 0\ 0\ 0\ 0\ 0\ 0\ 0\ 0\ 0\ 0\ 0\ 1\ 2\ 3\ 4\ 5\ 6\ 7\ 4\ 5\ 6\ 7\ 1\ 2\ 3\ 4\ 5\ 6\ 7\ 1\ 2\ 3\ 1\ 2\ 3\ 4\ 5\ 6\ 7\ 4\ 5\ 6\ 7\ 1\ 2\ 3 \\
3\ 4\ 5\ 6\ 7\ 1\ 2\ 1\ 2\ 3\ 4\ 5\ 6\ 7\ 0\ 0\ 0\ 0\ 0\ 0\ 0\ 0\ 0\ 0\ 0\ 0\ 0\ 0\ 1\ 2\ 3\ 4\ 5\ 6\ 7\ 3\ 4\ 5\ 6\ 7\ 1\ 2\ 3\ 4\ 5\ 6\ 7\ 1\ 2 \\
4\ 5\ 6\ 7\ 1\ 2\ 3\ 4\ 5\ 6\ 7\ 1\ 2\ 3\ 3\ 4\ 5\ 6\ 7\ 1\ 2\ 3\ 4\ 5\ 6\ 7\ 1\ 2\ 0\ 0\ 0\ 0\ 0\ 0\ 0\ 0\ 0\ 0\ 0\ 0\ 0\ 0\ 1\ 2\ 3\ 4\ 5\ 6\ 7
\end{array}$$

EXAMPLE 10. A-optimal and MV-optimal design for $p = 5$, $t = 4$ and $n = 48$:

$$\begin{array}{c}
0\ 0\ 0\ 0\ 1\ 2\ 3\ 4\ 1\ 2\ 3\ 4\ 1\ 2\ 3\ 4\ 0\ 0\ 0\ 0\ 1\ 2\ 3\ 4\ 1\ 2\ 3\ 4\ 1\ 2\ 3\ 4\ 1\ 2\ 3\ 4\ 0\ 0\ 0\ 0\ 1\ 2\ 3\ 4\ 1\ 2\ 3\ 4 \\
1\ 2\ 3\ 4\ 0\ 0\ 0\ 0\ 2\ 4\ 1\ 3\ 2\ 4\ 1\ 3\ 1\ 2\ 3\ 4\ 0\ 0\ 0\ 0\ 2\ 4\ 1\ 3\ 2\ 4\ 1\ 3\ 2\ 4\ 1\ 3\ 2\ 4\ 1\ 3\ 0\ 0\ 0\ 0\ 2\ 4\ 1\ 3 \\
2\ 4\ 1\ 3\ 2\ 4\ 1\ 3\ 0\ 0\ 0\ 0\ 3\ 1\ 4\ 2\ 2\ 4\ 1\ 3\ 2\ 4\ 1\ 3\ 0\ 0\ 0\ 0\ 3\ 1\ 4\ 2\ 3\ 1\ 4\ 2\ 3\ 1\ 4\ 2\ 3\ 1\ 4\ 2\ 0\ 0\ 0\ 0 \\
3\ 1\ 4\ 2\ 3\ 1\ 4\ 2\ 3\ 1\ 4\ 2\ 0\ 0\ 0\ 0\ 3\ 1\ 4\ 2\ 3\ 1\ 4\ 2\ 3\ 1\ 4\ 2\ 0\ 0\ 0\ 0\ 0\ 0\ 0\ 0\ 4\ 3\ 2\ 1\ 4\ 3\ 2\ 1\ 4\ 3\ 2\ 1 \\
4\ 3\ 2\ 1\ 4\ 3\ 2\ 1\ 4\ 3\ 2\ 1\ 4\ 3\ 2\ 1\ 4\ 3\ 2\ 1\ 4\ 3\ 2\ 1\ 4\ 3\ 2\ 1\ 4\ 3\ 2\ 1\ 4\ 3\ 2\ 1\ 0\ 0\ 0\ 0\ 0\ 0\ 0\ 0\ 0\ 0\ 0\ 0
\end{array}$$

**6. Discussion and closing remarks.** In this article we imposed two conditions on the class of competing designs: (i) In each design the control treatment should appear equally often in all $p$ periods. (ii) In each design no treatment is allowed to be immediately followed by itself in any experimental unit (i.e., $m_{ii} = 0$, $i = 0, 1, \dots, t$). Although the optimality of the newly discovered designs is over this restricted class, we are highly confident that these designs are very efficient over the entire class of designs. Indeed, for some $t$, $p$ and $n$ these designs could be optimal over the entire class. To support this optimism, we can see from Lemma 3 that for any design $d$ in which condition (i) is not satisfied, (3.2) may become a strict inequality. It is possible that this condition is a necessary condition for the optimal design in the entire class and therefore the restricted condition (i) could be removed. As for condition (ii), if $m_{dii} > 0$ for some $1 \leq i \leq t$, $x_0$ in Lemma 4 may not yield a significantly bigger value since the gain in the second part will be reduced by the loss in the first part. Consequently, we think the restricted condition $m_{dii} = 0$ for all $1 \leq i \leq t$ could indeed be a necessary condition for the optimal design, or at least the gain with the condition removed is very little. Further, if $m_{d00} > 0$, the numerator of $\Theta_2$ becomes $n(p-1)(\sum n_{d0u}\tilde{n}_{d0u} - pm_{d00})$ and $\Theta_2$ could be zero or negative. Thus, the corresponding inequalities (4.3) and (4.4) may not be greater than zero for a general design $d$ and consequently $\mathrm{Tr}(M_d^{-1})$ may not be minimized when $\sum n_{d0u}\tilde{n}_{d0u}$ or $\sum \tilde{n}_{d0u}^2$ is minimized for fixed $r_{d0}$. In this situation, we may have to consider the relationship among $\sum n_{d0u}^2$, $\sum n_{d0u}\tilde{n}_{d0u}$ and $\sum \tilde{n}_{d0u}^2$. We may find optimal designs for some special parameters, but we feel it will be very difficult to find optimal designs for the general parameters $p \leq t+1$.

It is natural to postulate that there could be better designs than those we have identified in this paper. While we do not have any general evidence for that, what so far we can say is this. The designs characterized in this paper are highly efficient if not optimal. We have substantial numerical evidence in support of this. Here is a typical example from the assorted examples that we have produced in our ongoing research in this area. Consider the case $t = 7$, $p = 4$ and $n = 28$. We used a computer along with some algebraic methods and searched for the possible lower bound of $\mathrm{Tr}(M_d^{-1})$ under the subclass of designs in which the control treatment appears equally often in $p$ periods. We found that there could be a design with $\mathrm{Tr}(M_d^{-1})$ equal to 1.02252. While we are not sure if there is a design with such a trace, let us assume that there is one. Clearly this hypothetical design is better than the A-optimal design $d^*$ with $\mathrm{Tr}(M_d^{-1}) = 1.02327$ which we displayed in Example 5 within the subclass $\Lambda_{t+1,n,p}$ of designs for these $t$, $p$ and $n$. But note that $d^*$ in Example 5 is at least 99.9% efficient in the larger class without the restriction $m_{ii} = 0$, $i = 0, 1, \dots, t$. This is not an isolated case, and as we mentioned, we have observed this phenomenon in many cases.

Another important issue worth discussing here is the status of the model robustness of the optimal and efficient designs discovered in this article. We used carryover Model (2.1) while searching for optimal and efficient designs. It is quite possible that upon data collection and analysis we may discover that another model could be more appropriate than Model (2.1). Consequently, the optimal design which was used under the postulated Model (2.1) might no longer be optimal or even efficient under the model specified in light of the data. Thus it is prudent to recommend a design to the experimenter which is optimal or at least efficient under several different likely models. For crossover studies several simple and lower models than Model (2.1) have been used in practice. We shall explore here how our designs in this article perform under the model without carryover effects (two-way elimination model), the model which contains only direct treatment effects and subject effects (one-way elimination model) and the model which contains only direct treatment effects (zero-elimination model). In Table 1 we have listed 15 designs of which 10 have already been displayed in Section 5 and the remaining 5 can be obtained by the procedures in this paper. This table lists the efficiency of each design under the carryover Model (2.1) in the subclass $\Lambda_{t+1,n,p}$ of designs ($e_c$), as well as zero-way ($e_0$), one-way ($e_1$), and two-way ($e_2$) elimination models in the unrestricted class of designs. The efficiencies $e_0$, $e_1$ and $e_2$ are based on the result of Hedayat, Jacroux and Majumdar (1988). It is clear that for most of these designs, if the design is A-optimal under the carryover model, then it is also A-optimal under one-way and two-way elimination models, and highly efficient under the zero-way elimination model.

TABLE 1
*Efficiencies under different models*

| **Design** | $p$ | $t$ | $n$ | $r_0$ | $e_c$ (%) | $e_0$ (%) | $e_1$ (%) | $e_2$ (%) |
|---|---|---|---|---|---|---|---|---|
| 1 | 3 | 2 | 6 | 6 | 100 | 97.50 | 100 | 100 |
| 2 | 3 | 3 | 9 | 9 | 100 | 100 | 100 | 100 |
| 3 | 3 | 4 | 36 | 36 | 100 | 100 | 100 | 100 |
| 4 | 3 | 5 | 30 | 30 | 100 | 99.84 | 99.85 | 99.85 |
| 5 | 3 | 7 | 49 | 42 | 100 | 99.98 | 99.97 | 99.97 |
| 6 | 4 | 3 | 4 | 4 | 100 | 94.4 | 98.75 | 98.75 |
| 7 | 4 | 4 | 16 | 16 | 100 | 96.55 | 100 | 100 |
| 8 | 4 | 5 | 40 | 40 | 100 | 98.18 | 100 | 100 |
| 9 | 4 | 6 | 40 | 40 | 100 | 99.16 | 100 | 100 |
| 10 | 4 | 7 | 28 | 28 | 100 | 99.82 | 100 | 100 |
| 11 | 4 | 9 | 48 | 48 | 100 | 100 | 100 | 100 |
| 12 | 5 | 4 | 48 | 60 | 100 | 96.43 | 97.56 | 97.56 |
| 13 | 5 | 5 | 50 | 50 | 99.50 | 93.10 | 96.87 | 96.87 |
| 14 | 5 | 6 | 30 | 30 | 99.30 | 95.23 | 98.53 | 98.53 |
| 15 | 5 | 7 | 70 | 70 | 100 | 96.68 | 99.47 | 99.47 |

Even an efficient design under the carryover model remains highly efficient under the other three models. Notice that while Design 13 is highly efficient under the carryover model, it is not that efficient under other models. The reason is this: For computational simplicity this design was constructed under $r_{d0} = 50$ while our computation showed that $r_{d^*0}$ is 60. Therefore, it is very safe to conclude that optimal and efficient crossover designs which are constructed in this article remain optimal or highly efficient under lower case models which are discussed here. Our research effort in this area is continuing.

## APPENDIX

PROPOSITION A.1. *For any $d \in \Lambda_{t+1,n,p}$, we have*

$$\tilde{r}_{d0} \le \frac{n(p-1)}{2}.$$

PROOF. Since the control treatment appears equally often in periods, we have $\tilde{r}_{d0} = (p-1)l_{d01}$. Also $l_{d01} \le n/2$ since no treatment (either test treatment or control) is followed by itself. Thus we obtain the conclusion. □

PROPOSITION A.2. *For any $d \in \Lambda_{t+1,n,p}$, where $p \le t+1$, we have*

$$\sum_{u=1}^{n} n_{d0u}\tilde{n}_{d0u} \le t[n(p-1) - \tilde{r}_{d0}].$$

PROOF. For given $\tilde{r}_{d0}$, since $\sum_{u=1}^{n} \tilde{n}_{d0u} = \tilde{r}_{d0}$ and $0 \le \tilde{n}_{d0u} \le \frac{p}{2}$, we have

$$\sum_{u=1}^{n} \tilde{n}_{d0u}^2 \le \begin{cases} \dfrac{p\tilde{r}_{d0}}{2}, & \text{when } p \ge 4, \\ \tilde{r}_{d0}, & \text{when } p = 3. \end{cases} \qquad \text{(A.1)}$$

On the other hand, we have $\sum_{u=1}^{n} n_{d0u}\tilde{n}_{d0u} \le \sum_{u=1}^{n} \tilde{n}_{d0u}(\tilde{n}_{d0u}+1)$. Notice that $3 \le p \le t+1$; then by applying (A.1) and Proposition A.1, we obtain the conclusion. □

In the next three propositions, $\xi_1$, $\xi_2$, $\xi_3$, $\widetilde{\Delta}_1$, $\widetilde{\Delta}_2$, $\widetilde{\Theta}_1$ and $\widetilde{\Theta}_2$ have the same definitions as those in the proof of Lemma 5. Notice that for any $d \in \Lambda_{t+1,n,p}$ and given $r_{d0}$, we have

$$\xi_1 \ge \min_d \left( \sum_{u=1}^{n} n_{d0u}^2 \right) \ge r_{d0},$$

$$nt(p-1) - t\tilde{r}_{d0} \ge \xi_2 \ge \min_d \left( \sum_{u=1}^{n} n_{d0u}\tilde{n}_{d0u} \right) \ge \tilde{r}_{d0},$$

$$\xi_3 \geq \min_d \left( \sum_{u=1}^{n} \tilde{n}_{d0u}^2 \right) \geq \tilde{r}_{d0}.$$

PROPOSITION A.3. *For any $d \in \Lambda_{t+1,n,p}$, where $3 \leq p \leq t+1$, we have*

$$\widetilde{\Delta}_1 \geq (t-1)\widetilde{\Delta}_2. \tag{A.2}$$

PROOF. Since $nt(p-1) - t\tilde{r}_{d0} \geq \xi_2 \geq \tilde{r}_{d0}$ and $\xi_3 \geq \tilde{r}_{d0}$, so we have

$$\begin{aligned} &\frac{(nt(p-1) - t\tilde{r}_{d0} - \xi_2)^2}{n(p-1)(pt-t-1) - (pt-t+p-2)\tilde{r}_{d0} + \xi_3} \\ &\qquad \leq \frac{(nt(p-1) - (t+1)\tilde{r}_{d0})^2}{n(p-1)(pt-t-1) - (pt-t+p-3)\tilde{r}_{d0}} \end{aligned}$$

and

$$\frac{n(p-1)(t-1)(\xi_2)^2}{np(p-1)\tilde{r}_{d0} - \tilde{r}_{d0}^2 - n(p-1)\xi_3} \geq \frac{n(p-1)(t-1)\tilde{r}_{d0}}{n(p-1)^2 - \tilde{r}_{d0}}.$$

We also notice that $tp(r_{d0} - \frac{1}{p}\xi_1) \leq t(pr_{d0} - \frac{r_{d0}^2}{n})$. By direct calculation it is sufficient to show

$$\begin{aligned} &t\left(pr_{d0} - \frac{r_{d0}^2}{n}\right) - t(p-1)(np - r_{d0}) - \frac{n(p-1)(t-1)\tilde{r}_{d0}}{n(p-1)^2 - \tilde{r}_{d0}} \\ &\qquad + \frac{(nt(p-1) - (t+1)\tilde{r}_{d0})^2}{n(p-1)(pt-t-1) - (pt-t+p-3)\tilde{r}_{d0}} \leq 0. \end{aligned} \tag{A.3}$$

The left-hand side of the preceding expression can be written as

$$\begin{aligned} -n\biggl(t\left(p - \frac{r_{d0}}{n}\right)&\left(p - 1 - \frac{r_{d0}}{n}\right) + \frac{(p-1)(t-1)(\tilde{r}_{d0}/n)}{(p-1)^2 - (\tilde{r}_{d0}/n)} \\ &- \frac{(t(p-1) - (t+1)(\tilde{r}_{d0}/n))^2}{(p-1)(pt-t-1) - (pt-t+p-3)(\tilde{r}_{d0}/n)}\biggr). \end{aligned}$$

From the proof of Proposition A.1, we know that $\frac{r_{d0}}{n} \leq \frac{\tilde{r}_{d0}}{n} + 1/2$. Define $x = \frac{\tilde{r}_{d0}}{n}$; then $\frac{r_{d0}}{n} \leq x + 1/2$. Simple algebra can show that (A.3) is equivalent to $f(x) > 0$ when $0 \leq x \leq (p-1)/2$, where

$$\begin{aligned} f(x) = t(2p-1-2x)(2p-3-2x) &+ \frac{4(t-1)x}{p-1} \\ &- \frac{4(t(p-1) - (t+1)x)^2}{(p-1)(pt-t-1) - (pt-t+p-3)x}. \end{aligned}$$

It can be shown that $(p-1)(pt-t-1)-(pt-t+p-3)x>0$ when $0\le x\le (p-1)/2$, so it is equivalent to showing $g(x)>0$ when $0\le x\le (p-1)/2$. Here

$$g(x)=\left[t(2p-1-2x)(2p-3-2x)+\frac{4(t-1)x}{p-1}\right]\\ \times[(p-1)(pt-t-1)-(pt-t+p-3)x]-4(t(p-1)-(t+1)x)^2.$$

It can be checked that $g(0)>0$ and $g(\frac{p-1}{2})>0$. If we can show that $g(x)$ is a monotone function when $0\le x\le (p-1)/2$, then we have reached the conclusion. By direct calculation, we have

$$\begin{aligned} g'(x)=&-12t(pt-t+p-3)x^2\\ &+\left[t^2(24p^2-48p+8)+t(16p^2-72p+40)+\frac{16(t-1)}{p-1}\right]x\\ &-t^2(12p^3-36p^2+23p+1)-t(4p^3-28p^2+39p-9)+4.\end{aligned}$$

Notice that $g'(x)$ is an increasing function when $x\in(-\infty,M)$, where

$$M=\frac{t^2(24p^2-48p+8)+t(16p^2-72p+40)+16(t-1)/(p-1)}{24t(pt-t+p-3)}.$$

It can be verified that $M>(p-1)/2$. So $g'(x)$ is an increasing function when $x\in[0,(p-1)/2]$. Also we can verify that

$$g'\left(\frac{p-1}{2}\right)=t^2(-3p^3+9p^2-4p-2)+t(p^3-p^2-4p+6)-4.$$

Notice that $-3p^3+9p^2-4p-2<0$ when $p\ge 3$; thus

$$\begin{aligned} g'\left(\frac{p-1}{2}\right)&\le 2t(-3p^3+9p^2-4p-2)+t(p^3-p^2-4p+6)-4\\ &=t(-5p^3+17p^2-12p+2)-4\\ &<0.\end{aligned}$$

So we have $g'(x)<0$ when $x\in[0,(p-1)/2]$; thus $g(x)$ is a monotone function when $x\in[0,(p-1)/2]$. $\square$

PROPOSITION A.4. *For any $d\in\Lambda_{t+1,n,p}$, with $p\le t+1$, and $\frac{\tilde r_{d0}}{n}\in[0,\frac{p-1}{t+1}]$, we have*

$$\frac{\widetilde\Delta_1}{(t-1)\widetilde\Delta_2}\ge\frac{t(p-1)}{t(p-1)-1}. \tag{A.4}$$

PROOF. Since $nt(p-1) - t\tilde{r}_{d0} \geq \xi_2$, $\widetilde{\Delta}_1$ will be minimized when $\xi_1$, $\xi_2$ and $\xi_3$ are minimized, that is,

$$\widetilde{\Delta}_1 \geq t(p-1)(np - r_{d0}) - (p-1)r_{d0} - \frac{[nt(p-1) - (t+1)\tilde{r}_{d0}]^2}{n(p-1)(pt-t-1) - (pt-t+p-3)\tilde{r}_{d0}}. \tag{A.5}$$

Also we can see that $\widetilde{\Delta}_2$ will be maximized when $\xi_1$, $\xi_2$ and $\xi_3$ are minimized, that is,

$$\widetilde{\Delta}_2 \leq (p-1)r_{d0} - \frac{n(p-1)\tilde{r}_{d0}}{n(p-1)^2 - \tilde{r}_{d0}}. \tag{A.6}$$

It suffices to show that

$$\frac{1}{n}[(t(p-1)-1)\widetilde{\Delta}_1 - t(t-1)(p-1)\widetilde{\Delta}_2] \geq 0.$$

Applying (A.5) and (A.6), and noticing that $r_{d0} = \frac{p}{p-1}\tilde{r}_{d0}$, we can show that

$$\begin{aligned} &\frac{1}{n}[(t(p-1)-1)\widetilde{\Delta}_1 - t(t-1)(p-1)\widetilde{\Delta}_2] \\ &\quad \geq [t(p-1)-1]\biggl[pt(p-1) - p(t+1)\frac{\tilde{r}_{d0}}{n} \\ &\qquad - \frac{t^2(p-1)^2}{(p-1)(pt-t-1) - (pt-t+p-3)(\tilde{r}_{d0}/n)}\biggr] \\ &\qquad - tp(t-1)(p-1)\frac{\tilde{r}_{d0}}{n}. \end{aligned} \tag{A.7}$$

Since $0 < \frac{\tilde{r}_{d0}}{n} < \frac{p-1}{t+1}$, and the right-hand side of (A.7) is a decreasing function of $\tilde{r}_{d0}$, we further have

$$\begin{aligned} &\frac{1}{n}[(t(p-1)-1)\widetilde{\Delta}_1 - t(t-1)(p-1)\widetilde{\Delta}_2] \\ &\quad \geq [t(p-1)-1] \\ &\qquad \times \biggl[pt(p-1) - p(p-1) \\ &\qquad - \frac{t^2(p-1)^2}{(p-1)(pt-t-1) - (pt-t+p-3)(p-1)/(t+1)}\biggr] \\ &\qquad - tp(t-1)(p-1)\frac{p-1}{t+1}. \end{aligned} \tag{A.8}$$

When $p \geq 4$, by applying the condition that $p \leq t+1$, we can show that $(pt-t-1) - \frac{pt-t+p-3}{t+1} \geq \frac{t(p-1)}{3}$, $t(p-1)-1 > \frac{3t(p-1)}{4}$ and $p(p-1)(t-1) - 3t >$

$\frac{p(p-1)(t-1)}{3}$, so by (A.8) we have

$$\begin{aligned}
&\frac{1}{n}[(t(p-1)-1)\widetilde{\Delta}_1 - t(t-1)(p-1)\widetilde{\Delta}_2] \\
&\qquad \geq [t(p-1)-1][p(t-1)(p-1)-3t] - tp(t-1)(p-1)\frac{p-1}{t+1} \\
&\qquad \geq tp(p-1)^2(t-1)\left(\frac{1}{4} - \frac{1}{t+1}\right) \\
&\qquad \geq 0.
\end{aligned}$$

When $p=3$ and $t \geq 3$, notice that $1 + \frac{2t}{t+1} < t$. By (A.8) we have

$$\begin{aligned}
&\frac{1}{n}[(t(p-1)-1)\widetilde{\Delta}_1 - t(t-1)(p-1)\widetilde{\Delta}_2] \\
&\qquad > 6(2t-1)(t-1) - (2t-1)\frac{2t^2}{2t-t} - 12\frac{t(t-1)}{t+1} \\
&\qquad \geq (2t-1)(4t-6) - 3t(t-1) \\
&\qquad > 0.
\end{aligned}$$

When $p=3$ and $t=2$, due to (A.5) and (A.6), we have $\widetilde{\Delta}_1 \geq 12n - 6r_{d0} - \frac{(4n-3\tilde{r}_{d0})^2}{6n-4\tilde{r}_{d0}}$ and $\widetilde{\Delta}_2 \leq 2r_{d0} - \frac{2n\tilde{r}_{d0}}{4n-\tilde{r}_{d0}}$. Thus we have

$$\begin{aligned}
&\frac{1}{n}[(t(p-1)-1)\widetilde{\Delta}_1 - t(t-1)(p-1)\widetilde{\Delta}_2] \\
&\qquad \geq 36 - 26\frac{r_{d0}}{n} - 3\frac{(4-3\tilde{r}_{d0}/n)^2}{6-4\tilde{r}_{d0}/n} + \frac{8\tilde{r}_{d0}/n}{4-\tilde{r}_{d0}/n} \\
&\qquad \geq 36 - 37\frac{\tilde{r}_{d0}}{n} - 3\frac{(4-3\tilde{r}_{d0}/n)^2}{6-4\tilde{r}_{d0}/n} \\
&\qquad > 0.
\end{aligned}$$

The last inequality can be easily verified when $0 \leq \frac{\tilde{r}_{d0}}{n} \leq \frac{2}{3}$. $\square$

PROPOSITION A.5. *For any $d \in \Lambda_{t+1,n,p}$, with $p \leq t+1$, we have*

$$0 \leq \frac{pt-t-1}{t(p-1)}\widetilde{\Theta}_1 \leq \widetilde{\Theta}_2. \tag{A.9}$$

*Furthermore, when $\frac{\tilde{r}_{d0}}{n} \geq \frac{p-1}{t+1}$, we have*

$$\widetilde{\Theta}_2 \geq \widetilde{\Theta}_1. \tag{A.10}$$

PROOF. It is easy to verify that $\widetilde{\Theta}_1 \geq 0$. We only focus on the remaining inequalities. By an argument similar to that in Proposition A.3, we have

$$\frac{nt(p-1) - t\tilde{r}_{d0} - \xi_2}{n(p-1)(pt-t-1) - (pt-t+p-2)\tilde{r}_{d0} + \xi_3} \leq \frac{nt(p-1) - (t+1)\tilde{r}_{d0}}{n(p-1)(pt-t-1) - (pt-t+p-3)\tilde{r}_{d0}}$$

and

$$\frac{n(p-1)\xi_2}{np(p-1)\tilde{r}_{d0} - \tilde{r}_{d0}^2 - n(p-1)\xi_3} \geq \frac{n(p-1)}{n(p-1)^2 - \tilde{r}_{d0}}.$$

As for (A.9), it is sufficient to show that

$$\text{(A.11)} \qquad \frac{p-1}{(p-1)^2 - \tilde{r}_{d0}/n} \frac{(p-1)(pt-t-1) - (pt-t+p-3)\tilde{r}_{d0}/n}{t(p-1) - (t+1)\tilde{r}_{d0}/n} \geq \frac{t(p-1)-1}{t(p-1)}.$$

Direct calculations show that (A.11) is equivalent to

$$(t+1)[t(p-1)-1]\left(\frac{\tilde{r}_{d0}}{n}\right)^2 - [(p-1)^2(t^2+t-1) - t(p-1)]\frac{\tilde{r}_{d0}}{n} \leq 0,$$

which holds when $0 \leq \frac{\tilde{r}_{d0}}{n} \leq \frac{(p-1)[(p-1)(t^2+t-1)-t]}{(t+1)[t(p-1)-1]}$. Notice that $(p-1)(t^2+t-1) - t \geq t^2(p-1) - t$, so $\frac{(p-1)[(p-1)(t^2+t-1)-t]}{(t+1)[t(p-1)-1]} \geq \frac{(p-1)t}{t+1}$. From Proposition A.1, we know that $\frac{\tilde{r}_{d0}}{n} \leq \frac{p-1}{2}$; thus (A.9) holds.

For (A.10) it suffices to show that

$$\text{(A.12)} \qquad \frac{p-1}{(p-1)^2 - \tilde{r}_{d0}/n} \geq \frac{t(p-1) - (t+1)\tilde{r}_{d0}/n}{(p-1)(pt-t-1) - (pt-t+p-3)\tilde{r}_{d0}/n}.$$

Direct calculations show that (A.12) is equivalent to

$$(t+1)\left(\frac{\tilde{r}_{d0}}{n}\right)^2 - (t+2)(p-1)\frac{\tilde{r}_{d0}}{n} + (p-1)^2 \leq 0,$$

which holds when $\frac{\tilde{r}_{d0}}{n} \in [\frac{p-1}{t+1}, p-1]$. From Proposition A.1 we know that $\frac{\tilde{r}_{d0}}{n} \leq \frac{p-1}{2}$. Thus (A.10) holds when $\frac{\tilde{r}_{d0}}{n} \geq \frac{p-1}{t+1}$. □

PROPOSITION A.6. *For any design $d \in \Lambda_{t+1,n,p}$ and given $r_{d0}$,*

$$\text{(A.13)} \qquad \min_d \left( \sum_{u=1}^n n_{d0u}^2 \right) = r_{d0} + (2r_{d0} - n)\left[\frac{r_{d0}}{n}\right] - n\left[\frac{r_{d0}}{n}\right]^2,$$

$$\min_d\left(\sum_{u=1}^n n_{d0u}\tilde n_{d0u}\right) = \begin{cases} \tilde r_{d0} + \left(2\tilde r_{d0} - n + \dfrac{\tilde r_{d0}}{p-1}\right)\left[\dfrac{\tilde r_{d0}}{n}\right] - n\left[\dfrac{\tilde r_{d0}}{n}\right]^2, \\ \qquad \textit{when } \dfrac{\tilde r_{d0}}{p-1} < n - \tilde r_{d0} + n\left[\dfrac{\tilde r_{d0}}{n}\right], \\ 2\tilde r_{d0} + \dfrac{\tilde r_{d0}}{p-1} - n + \left(2\tilde r_{d0} - 2n + \dfrac{\tilde r_{d0}}{p-1}\right)\left[\dfrac{\tilde r_{d0}}{n}\right] - n\left[\dfrac{\tilde r_{d0}}{n}\right]^2, \\ \qquad \textit{otherwise}, \end{cases} \tag{A.14}$$

*and*

$$\min_d\left(\sum_{u=1}^n \tilde n_{d0u}^2\right) = \tilde r_{d0} + (2\tilde r_{d0} - n)\left[\frac{\tilde r_{d0}}{n}\right] - n\left[\frac{\tilde r_{d0}}{n}\right]^2. \tag{A.15}$$

*$\sum_{u=1}^n n_{d0u}^2$, $\sum_{u=1}^n n_{d0u}\tilde n_{d0u}$ and $\sum_{u=1}^n \tilde n_{d0u}^2$ can achieve their minimum values when $d$ is a totally balanced test-control incomplete crossover design.*

PROOF. It is straightforward to show (A.13) and (A.15). Here we will prove (A.14) only. First we notice that

$$\sum_{u=1}^n n_{d0u}\tilde n_{d0u} = \sum_{u=1}^n \tilde n_{d0u}^2 + \sum_{u\in\Gamma} \tilde n_{d0u}.$$

Here $\Gamma$ is the set of units which receives the control treatment in the last period. So there are $r_{d0} - \tilde r_{d0} = \frac{\tilde r_{d0}}{p-1}$ units in $\Gamma$.

For any $\tilde n_{d0u}$, $u = 1, \ldots, n$, $\sum_{u=1}^n n_{d0u}\tilde n_{d0u}$ will be minimized when we put the $\frac{\tilde r_{d0}}{p-1}$ smallest values among all $\tilde n_{d0u}$ into $\Gamma$. Next, we will show that for a given $\tilde r_{d0}$, $\sum_{u=1}^n n_{d0u}\tilde n_{d0u}$ will be minimized when $d$ is a balanced test-control incomplete block design for the direct effects and carryover effects.

Suppose that there are some $\tilde n_{d0u}$'s, say $\tilde n_{d01}$ and $\tilde n_{d02}$, such that $\tilde n_{d01} - \tilde n_{d02} \ge 2$. Then we can replace $\tilde n_{d01}$ by $\tilde n'_{d01} = \tilde n_{d01} - 1$ and $\tilde n_{d02}$ by $\tilde n'_{d02} = \tilde n_{d02} + 1$ and keep the others unchanged. Direct calculations show that the value of $\sum_{u=1}^n \tilde n_{d0u}^2$ is decreased by at least 2. Meanwhile, $\sum_{u\in\Gamma}\tilde n_{d0u}$ is increased by at most 1. So $\sum_{u=1}^n n_{d0u}\tilde n_{d0u}$ is decreased by at least 1. Thus $\sum_{u=1}^n n_{d0u}\tilde n_{d0u}$ will be minimized when $d$ is a balanced test-control incomplete block design for carryover effects and $\Gamma$ is the set of units which has the $\frac{\tilde r_{d0}}{p-1}$ smallest values among $\tilde n_{d0u}$. When $d$ is a balanced test-control incomplete block design for the direct effects and carryover effects, it satisfies this condition. By direct calculations, we obtain the conclusion. □

**Acknowledgments.** We received numerous constructive and helpful comments from the Editor, the Associate Editor and two referees. These comments and suggestions substantially changed and improved the article. We gratefully thank them.

Department of Mathematics, Statistics,
and Computer Science
University of Illinois at Chicago
851 South Morgan Street
Chicago, Illinois 60607-7045
USA
e-mail: hedayat@uic.edu

Department of Statistics
University of Nebraska-Lincoln
Lincoln, Nebraska 68583-0712
USA
e-mail: myang2@unl.edu